# Bäcklund Transformations: Some Old and New Perspectives

C. J. Papachristou *,   A. N. Magoulas **

* Department of Physical Sciences, Hellenic Naval Academy, Piraeus 18539, Greece
E-mail:  papachristou@snd.edu.gr

** Department of Electrical Engineering, Hellenic Naval Academy, Piraeus 18539, Greece
E-mail:  aris@snd.edu.gr

**Abstract**

Bäcklund transformations (BTs) are traditionally regarded as a tool for integrating nonlinear partial differential equations (PDEs). Their use has been recently extended, however, to problems such as the construction of recursion operators for symmetries of PDEs, as well as the solution of linear systems of PDEs. In this article, the concept and some applications of BTs are reviewed. As an example of an integrable linear system of PDEs, the Maxwell equations of electromagnetism are shown to constitute a BT connecting the wave equations for the electric and the magnetic field; plane-wave solutions of the Maxwell system are constructed in detail. The connection between BTs and recursion operators is also discussed.

## 1.  Introduction

*Bäcklund transformations* (BTs) were originally devised as a tool for obtaining solutions of nonlinear partial differential equations (PDEs) (see, e.g., [1] and the references therein). They were later also proven useful as *recursion operators* for constructing infinite sequences of nonlocal symmetries and conservation laws of certain PDEs [2–6].

In simple terms, a BT is a system of PDEs connecting two fields that are required to independently satisfy two respective PDEs [say, (*a*) and (*b*)] in order for the system to be integrable for either field. If a solution of PDE (*a*) is known, then a solution of PDE (*b*) is obtained simply by integrating the BT, without having to actually solve the latter PDE (which, presumably, would be a much harder task). In the case where the PDEs (*a*) and (*b*) are identical, the *auto-BT* produces new solutions of PDE (*a*) from old ones.

As described above, a BT is an auxiliary tool for finding solutions of a given (usually nonlinear) PDE, using known solutions of the same or another PDE. But, what if the BT itself is the differential system whose solutions we are looking for? As it turns out, to solve the problem we need to have parameter-dependent solutions of *both* PDEs (*a*) and (*b*) at hand. By properly matching the parameters (provided this is possible) a solution of the given system is obtained.

The above method is particularly effective in *linear* problems, given that parametric solutions of linear PDEs are generally not hard to find. An important paradigm of a BT associated with a linear problem is offered by the Maxwell system of equations of electromagnetism [7,8]. As is well known, the consistency of this system demands that both the electric and the magnetic field independently satisfy a respective wave equation. These equations have known, parameter-dependent solutions; namely, monochromatic plane waves with arbitrary amplitudes, frequencies and wave vectors (the "parameters" of the problem). By inserting these solutions into the Maxwell system, one may find the appropriate expressions for the "parameters" in order for the plane waves to also be solutions of Maxwell's equations; that is, in order to represent an actual electromagnetic field.

This article, written for educational purposes, is an introduction to the concept of a BT and its application to the solution of PDEs or systems of PDEs. Both "classical" and novel views of a BT are discussed, the former view predominantly concerning integration of nonlinear PDEs while the latter one being applicable mostly to linear systems of PDEs. The article is organized as follows:

In Section 2 we review the classical concept of a BT. The solution-generating process by using a BT is demonstrated in a number of examples.

In Sec. 3 a different perception of a BT is presented, according to which it is the BT itself whose solutions are sought. The concept of *conjugate solutions* is introduced.

As an example, in Secs. 4 and 5 the Maxwell equations in empty space and in a linear conducting medium, respectively, are shown to constitute a BT connecting the wave equations for the electric and the magnetic field. Following [7], the process of constructing plane-wave solutions of this BT is presented in detail. This process is, of course, a familiar problem of electrodynamics but is seen here under a new perspective by employing the concept of a BT.

Finally, in Sec. 6 we briefly review the connection between BTs and recursion operators for generating infinite sequences of nonlocal symmetries of PDEs.

## 2. Bäcklund Transformations: Classical Viewpoint

Consider two PDEs $P[u]=0$ and $Q[v]=0$ for the unknown functions $u$ and $v$, respectively. The expressions $P[u]$ and $Q[v]$ may contain the corresponding variables $u$ and $v$, as well as partial derivatives of $u$ and $v$ with respect to the independent variables. For simplicity, we assume that $u$ and $v$ are functions of only two variables $x$, $t$. Partial derivatives with respect to these variables will be denoted by using subscripts: $u_x$, $u_t$, $u_{xx}$, $u_{tt}$, $u_{xt}$, etc.

Independently, for the moment, also consider a pair of coupled PDEs for $u$ and $v$:

$$B_1[u,v] = 0 \quad (a) \qquad B_2[u,v] = 0 \quad (b) \tag{1}$$

where the expressions $B_i[u,v]$ ($i=1,2$) may contain $u$, $v$ as well as partial derivatives of $u$ and $v$ with respect to $x$ and $t$. We note that $u$ appears in both equations (*a*) and (*b*). The question then is: if we find an expression for $u$ by integrating (*a*) for a given $v$, will it match the corresponding expression for $u$ found by integrating (*b*) for the same $v$? The answer is that, in order that (*a*) and (*b*) be consistent with each other for solution for $u$, the function $v$ must be properly chosen so as to satisfy a certain *consistency condition* (or *integrability condition* or *compatibility condition*).

By a similar reasoning, in order that (*a*) and (*b*) in (1) be mutually consistent for solution for $v$, for some given $u$, the function $u$ must now itself satisfy a corresponding integrability condition.

If it happens that the two consistency conditions for integrability of the system (1) are precisely the PDEs $P[u]=0$ and $Q[v]=0$, we say that the above system constitutes a *Bäcklund transformation* (BT) connecting solutions of $P[u]=0$ with solutions of $Q[v]=0$. In the special case where $P\equiv Q$, i.e., when $u$ and $v$ satisfy *the same* PDE, the system (1) is called an *auto-Bäcklund* transformation (auto-BT) for this PDE.

Suppose now that we seek solutions of the PDE $P[u]=0$. Assume that we are able to find a BT connecting solutions $u$ of this equation with solutions $v$ of the PDE $Q[v]=0$ (if $P\equiv Q$, the auto-BT connects solutions $u$ and $v$ of the same PDE) and let $v=v_0(x,t)$ be some known solution of $Q[v]=0$. The BT is then a system of PDEs for the unknown $u$,

$$B_i[u,v_0]=0, \quad i=1,2 \qquad (2)$$

The system (2) is integrable for $u$, given that the function $v_0$ satisfies *a priori* the required integrability condition $Q[v]=0$. The solution $u$ then of the system satisfies the PDE $P[u]=0$. Thus a solution $u(x,t)$ of the latter PDE is found without actually solving the equation itself, simply by integrating the BT (2) with respect to $u$. Of course, this method will be useful provided that integrating the system (2) for $u$ is simpler than integrating the PDE $P[u]=0$ itself. If the transformation (2) is an auto-BT for the PDE $P[u]=0$, then, starting with a known solution $v_0(x,t)$ of this equation and integrating the system (2), we find another solution $u(x,t)$ of the same equation.

Let us see some examples of the use of a BT to generate solutions of a PDE:

1. The *Cauchy-Riemann relations* of Complex Analysis,

$$u_x = v_y \quad (a) \qquad u_y = -v_x \quad (b) \qquad (3)$$

(here, the variable $t$ has been renamed $y$) constitute an auto-BT for the *Laplace equation*,

$$P[w] \equiv w_{xx} + w_{yy} = 0 \qquad (4)$$

Let us explain this: Suppose we want to solve the system (3) for $u$, for a given choice of the function $v(x,y)$. To see if the PDEs (*a*) and (*b*) match for solution for $u$, we must compare them in some way. We thus differentiate (*a*) with respect to $y$ and (*b*) with respect to $x$, and equate the mixed derivatives of $u$. That is, we apply the integrability condition $(u_x)_y = (u_y)_x$. In this way we eliminate the variable $u$ and find the condition that must be obeyed by $v(x,y)$:

$$P[v] \equiv v_{xx} + v_{yy} = 0 \ .$$

Similarly, by using the integrability condition $(v_x)_y = (v_y)_x$ to eliminate $v$ from the system (3), we find the necessary condition in order that this system be integrable for $v$, for a given function $u(x,y)$:

$$P[u] \equiv u_{xx} + u_{yy} = 0 \ .$$

In conclusion, the integrability of system (3) with respect to either variable requires that the other variable must satisfy the Laplace equation (4).

Let now $v_0(x,y)$ be a known solution of the Laplace equation (4). Substituting $v=v_0$ in the system (3), we can integrate this system with respect to $u$. It is not hard to show (by eliminating $v_0$ from the system) that the solution $u$ will also satisfy the Laplace equation (4). As an example, by choosing the solution $v_0(x,y)=xy$, we find a new solution $u(x,y)= (x^2 - y^2)/2 + C$.

2. The *Liouville equation* is written

$$P[u] \equiv u_{xt} - e^u = 0 \quad \Leftrightarrow \quad u_{xt} = e^u \tag{5}$$

Due to its nonlinearity, this PDE is hard to integrate directly. A solution is thus sought by means of a BT. We consider an auxiliary function $v(x,t)$ and an associated PDE,

$$Q[v] \equiv v_{xt} = 0 \tag{6}$$

We also consider the system of first-order PDEs,

$$u_x + v_x = \sqrt{2}\, e^{(u-v)/2} \quad (a) \qquad u_t - v_t = \sqrt{2}\, e^{(u+v)/2} \quad (b) \tag{7}$$

Differentiating the PDE (*a*) with respect to *t* and the PDE (*b*) with respect to *x*, and eliminating $(u_t - v_t)$ and $(u_x + v_x)$ in the ensuing equations with the aid of (*a*) and (*b*), we find that $u$ and $v$ satisfy the PDEs (5) and (6), respectively. Thus, the system (7) is a BT connecting solutions of (5) and (6). Starting with the trivial solution $v=0$ of (6), and integrating the system

$$u_x = \sqrt{2}\, e^{u/2}, \quad u_t = \sqrt{2}\, e^{u/2},$$

we find a nontrivial solution of (5):

$$u(x,t) = -2\ln\left(C - \frac{x+t}{\sqrt{2}}\right).$$

3. The "*sine-Gordon*" *equation* has applications in various areas of Physics, e.g., in the study of crystalline solids, in the transmission of elastic waves, in magnetism, in elementary-particle models, etc. The equation (whose name is a pun on the related linear Klein-Gordon equation) is written

$$P[u] \equiv u_{xt} - \sin u = 0 \quad \Leftrightarrow \quad u_{xt} = \sin u \tag{8}$$

The following system of equations is an auto-BT for the nonlinear PDE (8):

$$\frac{1}{2}(u+v)_x = a \sin\left(\frac{u-v}{2}\right), \quad \frac{1}{2}(u-v)_t = \frac{1}{a}\sin\left(\frac{u+v}{2}\right) \tag{9}$$

where $a\ (\neq 0)$ is an arbitrary real constant. [Because of the presence of *a*, the system (9) is called a *parametric* BT.] When $u$ is a solution of (8) the BT (9) is integrable for $v$, which, in turn, also is a solution of (8): $P[v]=0$; and vice versa. Starting with the trivial solution $v=0$ of $v_{xt}= \sin v$, and integrating the system

$$u_x = 2a \sin \frac{u}{2} \ , \quad u_t = \frac{2}{a} \sin \frac{u}{2} \ ,$$

we obtain a new solution of (8):

$$u(x,t) = 4 \arctan \left\{ C \exp \left( ax + \frac{t}{a} \right) \right\} \ .$$

## 3. Conjugate Solutions and Another View of a BT

As presented in the previous section, a BT is an auxiliary device for constructing solutions of a (usually nonlinear) PDE from known solutions of the same or another PDE. The converse problem, where solutions of the differential system representing the BT itself are sought, is also of interest, however, and has been recently suggested [7,8] in connection with the Maxwell equations (see subsequent sections).

To be specific, assume that we need to integrate a given system of PDEs connecting two functions $u$ and $v$:

$$B_i[u,v] = 0 \ , \quad i = 1, 2 \tag{10}$$

Suppose that the integrability of the system for both functions requires that $u$ and $v$ separately satisfy the respective PDEs

$$P[u] = 0 \quad (a) \qquad Q[v] = 0 \quad (b) \tag{11}$$

That is, the system (10) is a BT connecting solutions of the PDEs (11). Assume, now, that these PDEs possess known (or, in any case, easy to find) *parameter-dependent solutions* of the form

$$u = f(x, y; \alpha, \beta, \ldots) \ , \quad v = g(x, y; \kappa, \lambda, \ldots) \tag{12}$$

where $\alpha$, $\beta$, $\kappa$, $\lambda$, etc., are (real or complex) parameters. If values of these parameters can be determined for which $u$ and $v$ jointly satisfy the system (10), we say that the solutions $u$ and $v$ of the PDEs (11a) and (11b), respectively, are *conjugate through the BT* (10) (or *BT-conjugate*, for short). By finding a pair of BT-conjugate solutions one thus automatically obtains a solution of the system (10).

Note that solutions of *both* integrability conditions $P[u]=0$ and $Q[v]=0$ must now be known in advance! From the practical point of view the method is thus most applicable in *linear* problems, since it is much easier to find parameter-dependent solutions of the PDEs (11) in this case.

Let us see an example: Going back to the Cauchy-Riemann relations (3), we try the following parametric solutions of the Laplace equation (4):

$$u(x, y) = \alpha (x^2 - y^2) + \beta x + \gamma y \ ,$$
$$v(x, y) = \kappa xy + \lambda x + \mu y \ .$$

Substituting these into the BT (3), we find that $\kappa=2\alpha$, $\mu=\beta$ and $\lambda=-\gamma$. Therefore, the solutions

$$u(x,y) = \alpha(x^2 - y^2) + \beta x + \gamma y ,$$
$$v(x,y) = 2\alpha xy - \gamma x + \beta y$$

of the Laplace equation are BT-conjugate through the Cauchy-Riemann relations.

As a counter-example, let us try a different combination:

$$u(x,y) = \alpha xy , \quad v(x,y) = \beta xy .$$

Inserting these into the system (3) and taking into account the independence of $x$ and $y$, we find that the only possible values of the parameters $\alpha$ and $\beta$ are $\alpha=\beta=0$, so that $u(x,y)=v(x,y)=0$. Thus, no non-trivial BT-conjugate solutions exist in this case.

## 4. Example: The Maxwell Equations in Empty Space

An example of an integrable linear system whose solutions are of physical interest is furnished by the *Maxwell equations* of electrodynamics. Interestingly, as noted recently [7], the Maxwell system has the property of a BT whose integrability conditions are the electromagnetic (e/m) wave equations that are separately valid for the electric and the magnetic field. These equations possess parameter-dependent solutions that, by a proper choice of the parameters, can be made BT-conjugate through the Maxwell system. In this and the following section we discuss the BT property of the Maxwell equations in vacuum and in a conducting medium, respectively.

In *empty space*, where no charges or currents (whether free or bound) exist, the Maxwell equations are written (in S.I. units) [9]

$$(a) \quad \vec{\nabla} \cdot \vec{E} = 0 \qquad (c) \quad \vec{\nabla} \times \vec{E} = -\frac{\partial \vec{B}}{\partial t}$$
$$(b) \quad \vec{\nabla} \cdot \vec{B} = 0 \qquad (d) \quad \vec{\nabla} \times \vec{B} = \varepsilon_0 \mu_0 \frac{\partial \vec{E}}{\partial t} \tag{13}$$

where $\vec{E}$ and $\vec{B}$ are the electric and the magnetic field, respectively. Here we have a system of four PDEs for two fields. The question is: what are the necessary conditions that each of these fields must satisfy in order for the system (13) to be self-consistent? In other words, what are the *consistency conditions* (or *integrability conditions*) for this system?

Guided by our experience from Sec. 2, to find these conditions we perform various differentiations of the equations of system (13) and require that certain differential identities be satisfied. Our aim is, of course, to eliminate one field (electric or magnetic) in favor of the other and find some higher-order PDE that the latter field must obey.

As can be checked, two differential identities are satisfied automatically in the system (13):

$$\vec{\nabla} \cdot (\vec{\nabla} \times \vec{E}) = 0 \ , \quad \vec{\nabla} \cdot (\vec{\nabla} \times \vec{B}) = 0 \ ,$$

$$(\vec{\nabla} \cdot \vec{E})_t = \vec{\nabla} \cdot \vec{E}_t \ , \quad (\vec{\nabla} \cdot \vec{B})_t = \vec{\nabla} \cdot \vec{B}_t \ .$$

Two others read

$$\vec{\nabla} \times (\vec{\nabla} \times \vec{E}) = \vec{\nabla}(\vec{\nabla} \cdot \vec{E}) - \nabla^2 \vec{E} \tag{14}$$

$$\vec{\nabla} \times (\vec{\nabla} \times \vec{B}) = \vec{\nabla}(\vec{\nabla} \cdot \vec{B}) - \nabla^2 \vec{B} \tag{15}$$

Taking the *rot* of (13*c*) and using (14), (13*a*) and (13*d*), we find

$$\nabla^2 \vec{E} - \varepsilon_0 \mu_0 \frac{\partial^2 \vec{E}}{\partial t^2} = 0 \tag{16}$$

Similarly, taking the *rot* of (13*d*) and using (15), (13*b*) and (13*c*), we get

$$\nabla^2 \vec{B} - \varepsilon_0 \mu_0 \frac{\partial^2 \vec{B}}{\partial t^2} = 0 \tag{17}$$

No new information is furnished by the remaining two integrability conditions,

$$(\vec{\nabla} \times \vec{E})_t = \vec{\nabla} \times \vec{E}_t \ , \quad (\vec{\nabla} \times \vec{B})_t = \vec{\nabla} \times \vec{B}_t \ .$$

Note that we have *uncoupled* the equations for the two fields in the system (13), deriving separate second-order PDEs for each field. Putting

$$\varepsilon_0 \mu_0 \equiv \frac{1}{c^2} \Leftrightarrow c = \frac{1}{\sqrt{\varepsilon_0 \mu_0}} \tag{18}$$

(where $c$ is the speed of light in vacuum) we rewrite (16) and (17) in wave-equation form:

$$\nabla^2 \vec{E} - \frac{1}{c^2} \frac{\partial^2 \vec{E}}{\partial t^2} = 0 \tag{19}$$

$$\nabla^2 \vec{B} - \frac{1}{c^2} \frac{\partial^2 \vec{B}}{\partial t^2} = 0 \tag{20}$$

We conclude that the Maxwell system (13) is a BT relating solutions of the e/m wave equations (19) and (20), these equations representing the integrability conditions of the BT. It should be noted that this BT is *not* an *auto*-BT! Indeed, although the PDEs (19) and (20) are of similar form, they concern *different* fields with different physical dimensions and physical properties.

The e/m wave equations admit plane-wave solutions of the form $\vec{F}(\vec{k}\cdot\vec{r}-\omega t)$, with

$$\frac{\omega}{k}=c \quad \text{where} \quad k=|\vec{k}| \tag{21}$$

The simplest such solutions are *monochromatic plane waves* of angular frequency $\omega$, propagating in the direction of the wave vector $\vec{k}$:

$$\begin{aligned}\vec{E}(\vec{r},t)&=\vec{E}_0\exp\{i(\vec{k}\cdot\vec{r}-\omega t)\} & (a)\\ \vec{B}(\vec{r},t)&=\vec{B}_0\exp\{i(\vec{k}\cdot\vec{r}-\omega t)\} & (b)\end{aligned} \tag{22}$$

where $\vec{E}_0$ and $\vec{B}_0$ are constant complex amplitudes. The constants appearing in the above equations (amplitudes, frequency and wave vector) can be chosen arbitrarily; thus they can be regarded as *parameters* on which the plane waves (22) depend.

We must note carefully that, although every pair of fields $(\vec{E},\vec{B})$ satisfying the Maxwell equations (13) also satisfies the wave equations (19) and (20), the converse is not true. Thus, the plane-wave solutions (22) are not *a priori* solutions of the Maxwell system (i.e., do not represent actual e/m fields). This problem can be taken care of, however, by a proper choice of the parameters in (22). To this end, we substitute the general solutions (22) into the BT (13) to find the extra conditions the latter system demands. By fixing the wave parameters, the two wave solutions in (22) will become *BT-conjugate* through the Maxwell system (13).

Substituting (22*a*) and (22*b*) into (13*a*) and (13*b*), respectively, and taking into account that $\vec{\nabla}e^{i\vec{k}\cdot\vec{r}}=i\vec{k}\,e^{i\vec{k}\cdot\vec{r}}$, we have

$$(\vec{E}_0\,e^{-i\omega t})\cdot\vec{\nabla}e^{i\vec{k}\cdot\vec{r}}=0 \;\Rightarrow\; (\vec{k}\cdot\vec{E}_0)\,e^{i(\vec{k}\cdot\vec{r}-\omega t)}=0,$$

$$(\vec{B}_0\,e^{-i\omega t})\cdot\vec{\nabla}e^{i\vec{k}\cdot\vec{r}}=0 \;\Rightarrow\; (\vec{k}\cdot\vec{B}_0)\,e^{i(\vec{k}\cdot\vec{r}-\omega t)}=0,$$

so that

$$\vec{k}\cdot\vec{E}_0=0, \quad \vec{k}\cdot\vec{B}_0=0. \tag{23}$$

Relations (23) reflect the fact that that the monochromatic plane e/m wave is a *transverse wave*.

Next, substituting (22*a*) and (22*b*) into (13*c*) and (13*d*), we find

$$e^{-i\omega t}(\vec{\nabla}e^{i\vec{k}\cdot\vec{r}})\times\vec{E}_0=i\omega\vec{B}_0\,e^{i(\vec{k}\cdot\vec{r}-\omega t)} \;\Rightarrow$$
$$(\vec{k}\times\vec{E}_0)\,e^{i(\vec{k}\cdot\vec{r}-\omega t)}=\omega\vec{B}_0\,e^{i(\vec{k}\cdot\vec{r}-\omega t)},$$

$$e^{-i\omega t}(\vec{\nabla}e^{i\vec{k}\cdot\vec{r}})\times\vec{B}_0=-i\omega\varepsilon_0\mu_0\vec{E}_0\,e^{i(\vec{k}\cdot\vec{r}-\omega t)} \;\Rightarrow$$
$$(\vec{k}\times\vec{B}_0)\,e^{i(\vec{k}\cdot\vec{r}-\omega t)}=-\frac{\omega}{c^2}\vec{E}_0\,e^{i(\vec{k}\cdot\vec{r}-\omega t)},$$

so that

$$\vec{k} \times \vec{E}_0 = \omega \vec{B}_0 \;, \quad \vec{k} \times \vec{B}_0 = -\frac{\omega}{c^2} \vec{E}_0 \tag{24}$$

We note that the fields $\vec{E}$ and $\vec{B}$ are normal to each other, as well as normal to the direction of propagation of the wave. We also remark that the two vector equations in (24) are not independent of each other, since, by cross-multiplying the first relation by $\vec{k}$, we get the second relation.

Introducing a unit vector $\hat{\tau}$ in the direction of the wave vector $\vec{k}$,

$$\hat{\tau} = \vec{k}/k \quad (k = |\vec{k}| = \omega/c) \;,$$

we rewrite the first of equations (24) as

$$\vec{B}_0 = \frac{k}{\omega} (\hat{\tau} \times \vec{E}_0) = \frac{1}{c} (\hat{\tau} \times \vec{E}_0) \;.$$

The BT-conjugate solutions in (22) are now written

$$\begin{aligned} \vec{E}(\vec{r},t) &= \vec{E}_0 \exp\{i(\vec{k} \cdot \vec{r} - \omega t)\} \;, \\ \vec{B}(\vec{r},t) &= \frac{1}{c} (\hat{\tau} \times \vec{E}_0) \exp\{i(\vec{k} \cdot \vec{r} - \omega t)\} = \frac{1}{c} \hat{\tau} \times \vec{E} \end{aligned} \tag{25}$$

As constructed, the complex vector fields in (25) satisfy the Maxwell system (13). Since this system is homogeneous linear with real coefficients, the real parts of the fields (25) also satisfy it. To find the expressions for the real solutions (which, after all, carry the physics of the situation) we take the simplest case of *linear polarization* and write

$$\vec{E}_0 = \vec{E}_{0,R} e^{i\alpha} \tag{26}$$

where the vector $\vec{E}_{0,R}$ as well as the number $\alpha$ are real. The *real* versions of the fields (25), then, read

$$\begin{aligned} \vec{E} &= \vec{E}_{0,R} \cos(\vec{k} \cdot \vec{r} - \omega t + \alpha) \;, \\ \vec{B} &= \frac{1}{c} (\hat{\tau} \times \vec{E}_{0,R}) \cos(\vec{k} \cdot \vec{r} - \omega t + \alpha) = \frac{1}{c} \hat{\tau} \times \vec{E} \end{aligned} \tag{27}$$

We note, in particular, that the fields $\vec{E}$ and $\vec{B}$ "oscillate" in phase.

Our results for the Maxwell equations in vacuum can be extended to the case of a *linear non-conducting medium* upon replacement of $\varepsilon_0$ and $\mu_0$ with $\varepsilon$ and $\mu$, respectively. The speed of propagation of the e/m wave is, in this case,

$$\upsilon = \frac{\omega}{k} = \frac{1}{\sqrt{\varepsilon \mu}} \;.$$

In the next section we study the more complex case of a linear medium having a finite conductivity.

## 5. Example: The Maxwell System for a Linear Conducting Medium

Consider a linear conducting medium of conductivity $\sigma$. In such a medium, Ohm's law is satisfied: $\vec{J}_f = \sigma \vec{E}$, where $\vec{J}_f$ is the free current density. The Maxwell equations take on the form [9]

$$
\begin{aligned}
&(a) \quad \vec{\nabla} \cdot \vec{E} = 0 \qquad (c) \quad \vec{\nabla} \times \vec{E} = -\frac{\partial \vec{B}}{\partial t} \\
&(b) \quad \vec{\nabla} \cdot \vec{B} = 0 \qquad (d) \quad \vec{\nabla} \times \vec{B} = \mu \sigma \vec{E} + \varepsilon \mu \frac{\partial \vec{E}}{\partial t}
\end{aligned}
\tag{28}
$$

By requiring satisfaction of the integrability conditions

$$
\vec{\nabla} \times (\vec{\nabla} \times \vec{E}) = \vec{\nabla}(\vec{\nabla} \cdot \vec{E}) - \nabla^2 \vec{E},
$$
$$
\vec{\nabla} \times (\vec{\nabla} \times \vec{B}) = \vec{\nabla}(\vec{\nabla} \cdot \vec{B}) - \nabla^2 \vec{B},
$$

we obtain the *modified wave equations*

$$
\begin{aligned}
\nabla^2 \vec{E} - \varepsilon \mu \frac{\partial^2 \vec{E}}{\partial t^2} - \mu \sigma \frac{\partial \vec{E}}{\partial t} = 0 \\
\nabla^2 \vec{B} - \varepsilon \mu \frac{\partial^2 \vec{B}}{\partial t^2} - \mu \sigma \frac{\partial \vec{B}}{\partial t} = 0
\end{aligned}
\tag{29}
$$

which must be separately satisfied by each field. As in Sec. 4, no further information is furnished by the remaining integrability conditions.

The linear differential system (28) is a BT relating solutions of the wave equations (29). As in the vacuum case, this BT is *not* an auto-BT. We now seek BT-conjugate solutions. As can be verified by direct substitution into equations (29), these PDEs admit parameter-dependent solutions of the form

$$
\begin{aligned}
\vec{E}(\vec{r},t) &= \vec{E}_0 \exp\{-s\,\hat{\tau}\cdot\vec{r} + i(\vec{k}\cdot\vec{r} - \omega t)\} \\
&= \vec{E}_0 \exp\left\{\left(i - \frac{s}{k}\right)\vec{k}\cdot\vec{r}\right\} \exp(-i\omega t), \\
\vec{B}(\vec{r},t) &= \vec{B}_0 \exp\{-s\,\hat{\tau}\cdot\vec{r} + i(\vec{k}\cdot\vec{r} - \omega t)\} \\
&= \vec{B}_0 \exp\left\{\left(i - \frac{s}{k}\right)\vec{k}\cdot\vec{r}\right\} \exp(-i\omega t)
\end{aligned}
\tag{30}
$$

where $\hat{\tau}$ is the unit vector in the direction of the wave vector $\vec{k}$:

$$\hat{\tau} = \vec{k}/k \quad (k = |\vec{k}| = \omega/v)$$

($v$ is the speed of propagation of the wave inside the conducting medium) and where, for given physical characteristics $\varepsilon$, $\mu$, $\sigma$ of the medium, the parameters $s$, $k$ and $\omega$ satisfy the algebraic system

$$s^2 - k^2 + \varepsilon\mu\omega^2 = 0 , \quad \mu\sigma\omega - 2sk = 0 \tag{31}$$

We note that, for arbitrary choices of the amplitudes $\vec{E}_0$ and $\vec{B}_0$, the vector fields (30) are not *a priori* solutions of the Maxwell system (28), thus are not BT-conjugate solutions. To obtain such solutions we substitute expressions (30) into the system (28). With the aid of the relation

$$\vec{\nabla} e^{\left(i - \frac{s}{k}\right)\vec{k}\cdot\vec{r}} = \left(i - \frac{s}{k}\right) \vec{k} \, e^{\left(i - \frac{s}{k}\right)\vec{k}\cdot\vec{r}}$$

one can show that (28*a*) and (28*b*) impose the conditions

$$\vec{k}\cdot\vec{E}_0 = 0 , \quad \vec{k}\cdot\vec{B}_0 = 0 \tag{32}$$

As in the vacuum case, the e/m wave in a conducting medium is a *transverse* wave.

By substituting (30) into (28*c*) and (28*d*), two more conditions are found:

$$(k + is)\hat{\tau}\times\vec{E}_0 = \omega\vec{B}_0 \tag{33}$$

$$(k + is)\hat{\tau}\times\vec{B}_0 = -(\varepsilon\mu\omega + i\mu\sigma)\vec{E}_0 \tag{34}$$

Note, however, that (34) is not an independent equation since it can be reproduced by cross-multiplying (33) by $\hat{\tau}$, taking into account the algebraic relations (31).

The BT-conjugate solutions of the wave equations (29) are now written

$$\begin{aligned}\vec{E}(\vec{r},t) &= \vec{E}_0 \, e^{-s\hat{\tau}\cdot\vec{r}} e^{i(\vec{k}\cdot\vec{r} - \omega t)} , \\ \vec{B}(\vec{r},t) &= \frac{k + is}{\omega} (\hat{\tau}\times\vec{E}_0) e^{-s\hat{\tau}\cdot\vec{r}} e^{i(\vec{k}\cdot\vec{r} - \omega t)}\end{aligned} \tag{35}$$

To find the corresponding real solutions, we assume linear polarization of the wave, as before, and set

$$\vec{E}_0 = \vec{E}_{0,R} \, e^{i\alpha} .$$

We also put

$$k + is = |k + is|e^{i\varphi} = \sqrt{k^2 + s^2} \, e^{i\varphi} ; \quad \tan\varphi = s/k .$$

Taking the real parts of equations (35), we finally have:

$$\vec{E}(\vec{r},t) = \vec{E}_{0,R}\, e^{-s\hat{\tau}\cdot\vec{r}} \cos(\vec{k}\cdot\vec{r} - \omega t + \alpha),$$

$$\vec{B}(\vec{r},t) = \frac{\sqrt{k^2 + s^2}}{\omega} (\hat{\tau} \times \vec{E}_{0,R})\, e^{-s\hat{\tau}\cdot\vec{r}} \cos(\vec{k}\cdot\vec{r} - \omega t + \alpha + \varphi).$$

As an exercise, the student may show that these results reduce to those for a linear non-conducting medium (cf. Sec. 4) in the limit $\sigma \to 0$.

## 6. BTs as Recursion Operators

The concept of symmetries of PDEs was discussed in [1]. Let us review the main facts:

Consider a PDE $F[u]=0$, where, for simplicity, $u=u(x,t)$. A transformation

$$u(x,t) \to u'(x,t)$$

from the function $u$ to a new function $u'$ represents a *symmetry* of the given PDE if the following condition is satisfied: $u'(x,t)$ is a solution of $F[u]=0$ if $u(x,t)$ is a solution. That is,

$$F[u']=0 \quad \text{when} \quad F[u]=0 \tag{36}$$

An *infinitesimal symmetry transformation* is written

$$u' = u + \delta u = u + \alpha Q[u] \tag{37}$$

where $\alpha$ is an infinitesimal parameter. The function $Q[u] \equiv Q(x, t, u, u_x, u_t, ...)$ is called the *symmetry characteristic* of the transformation (37).

In order that a function $Q[u]$ be a symmetry characteristic for the PDE $F[u]=0$, it must satisfy a certain PDE that expresses the *symmetry condition* for $F[u]=0$. We write, symbolically,

$$S(Q;u) = 0 \quad \text{when} \quad F[u]=0 \tag{38}$$

where the expression $S$ depends *linearly* on $Q$ and its partial derivatives. Thus, (38) is a linear PDE for $Q$, in which equation the variable $u$ enters as a sort of parametric function that is required to satisfy the PDE $F[u]=0$.

A *recursion operator* $\hat{R}$ [10] is a linear operator which, acting on a symmetry characteristic $Q$, produces a new symmetry characteristic $Q' = \hat{R}Q$. That is,

$$S(\hat{R}Q;u) = 0 \quad \text{when} \quad S(Q;u) = 0 \tag{39}$$

It is not too difficult to show that *any power of a recursion operator also is a recursion operator*. This means that, starting with any symmetry characteristic $Q$, one may in principle obtain an infinite set of characteristics (thus, an infinite number of symmetries) by repeated application of the recursion operator.

A new approach to recursion operators was suggested in the early 1990s [2,3] (see also [4-6]). According to this view, a recursion operator is an auto-BT for the linear PDE (38) expressing the symmetry condition of the problem; that is, a BT producing new solutions $Q'$ of (38) from old ones, $Q$. Typically, this type of BT produces *nonlocal* symmetries, i.e., symmetry characteristics depending on *integrals* (rather than derivatives) of $u$.

As an example, consider the *chiral field equation*

$$F[g] \equiv (g^{-1} g_x)_x + (g^{-1} g_t)_t = 0 \tag{40}$$

(as usual, subscripts denote partial differentiations) where $g$ is a *GL(n,C)*-valued function of $x$ and $t$ (i.e., an invertible complex $n \times n$ matrix, differentiable for all $x$, $t$).

Let $Q[g]$ be a symmetry characteristic of the PDE (40). It is convenient to put

$$Q[g] = g \Phi[g]$$

and write the corresponding infinitesimal symmetry transformation in the form

$$g' = g + \delta g = g + \alpha g \Phi[g] \tag{41}$$

The symmetry condition that $Q$ must satisfy will be a PDE linear in $Q$, thus in $\Phi$ also. As can be shown [4], this PDE is

$$S(\Phi; g) \equiv \Phi_{xx} + \Phi_{tt} + [g^{-1} g_x, \Phi_x] + [g^{-1} g_t, \Phi_t] = 0 \tag{42}$$

which must be valid when $F[g]=0$ (where, in general, $[A, B] \equiv AB - BA$ denotes the *commutator* of two matrices $A$ and $B$).

For a given $g$ satisfying $F[g]=0$, consider now the following system of PDEs for the matrix functions $\Phi$ and $\Phi'$:

$$\begin{aligned} \Phi'_x &= \Phi_t + [g^{-1} g_t, \Phi] \\ -\Phi'_t &= \Phi_x + [g^{-1} g_x, \Phi] \end{aligned} \tag{43}$$

The integrability condition $(\Phi'_x)_t = (\Phi'_t)_x$, together with the equation $F[g]=0$, require that $\Phi$ be a solution of (42): $S(\Phi; g) = 0$. Similarly, by the integrability condition $(\Phi_t)_x = (\Phi_x)_t$ one finds, after a lengthy calculation: $S(\Phi'; g) = 0$.

In conclusion, for any $g$ satisfying the PDE (40), the system (43) is a BT relating solutions $\Phi$ and $\Phi'$ of the symmetry condition (42) of this PDE; that is, relating different symmetries of the chiral field equation (40). Thus, if a symmetry characteristic $Q = g\Phi$ of (40) is known, a new characteristic $Q' = g\Phi'$ may be found by integrating the BT (43); the converse is also true. Since the BT (43) produces new symmetries from old ones, it may be regarded as a *recursion operator* for the PDE (40).

As an example, for any constant matrix *M* the choice Φ=*M* clearly satisfies the symmetry condition (42). This corresponds to the symmetry characteristic *Q*=*gM*. By integrating the BT (43) for Φ´, we get Φ´=[*X*, *M*] and *Q*´=*g*[*X*, *M*], where *X* is the "potential" of the PDE (40), defined by the system of PDEs

$$X_x = g^{-1} g_t \ , \quad -X_t = g^{-1} g_x \qquad (44)$$

Note the *nonlocal* character of the BT-produced symmetry *Q*´, due to the presence of the potential *X*. Indeed, as seen from (44), in order to find *X* one has to *integrate* the chiral field *g* with respect to the independent variables *x* and *t*. The above process can be continued indefinitely by repeated application of the recursion operator (43), leading to an infinite sequence of increasingly nonlocal symmetries.

## 7. Summary

Classically, Bäcklund transformations (BTs) have been developed as a useful tool for finding solutions of nonlinear PDEs, given that these equations are usually hard to solve by direct methods. By means of examples we saw that, starting with even the most trivial solution of a PDE, one may produce a highly nontrivial solution of this (or another) PDE by integrating the BT, without solving the original, nonlinear PDE directly (which, in most cases, is a much harder task).

A different use of BTs, that was recently proposed [7,8], concerns predominantly the solution of linear systems of PDEs. This method relies on the existence of parameter-dependent solutions of the linear PDEs expressing the integrability conditions of the BT. This time it is the BT itself (rather than its associated integrability conditions) whose solutions are sought.

An appropriate example for demonstrating this approach to the concept of a BT is furnished by the Maxwell equations of electromagnetism. We showed that this system of PDEs can be treated as a BT whose integrability conditions are the wave equations for the electric and the magnetic field. These wave equations have known, parameter-dependent solutions – monochromatic plane waves – with arbitrary amplitudes, frequencies and wave vectors playing the roles of the "parameters". By substituting these solutions into the BT, one may determine the required relations among the parameters in order that these plane waves also represent electromagnetic fields (i.e., in order that they be solutions of the Maxwell system). The results arrived at by this method are, of course, well known in advanced electrodynamics. The process of deriving them, however, is seen here in a new light by employing the concept of a BT.

BTs have also proven useful as *recursion operators* for deriving infinite sets of nonlocal symmetries and conservation laws of PDEs [2-6] (see also [11] and the references therein). Specifically, the BT produces an increasingly nonlocal sequence of symmetry characteristics, i.e., solutions of the linear equation expressing the symmetry condition (or "linearization") of a given PDE.

An interesting conclusion is that the concept of a BT, which has been proven useful for integrating nonlinear PDEs, may also have important applications in linear problems. Research on these matters is in progress.